\documentstyle[12pt,twoside]{article}
\newtheorem{theorem}{Theorem}[section]

\newtheorem{proposition}[theorem]{Proposition}

\newtheorem{lemma}[theorem]{Lemma}

\newcommand{\proof}{\noindent{\bf Proof. }}
\newcommand{\bbfR}{I\!\!R}
\renewcommand{\L}{{\cal L}}
\newcommand{\eLt}{{e^{-t\L}}}
\newcommand{\eLtt}{{e^{-(t-\tau)\L}}}
\newcommand{\e}{{\varepsilon}}

\newcommand{\Y}{{\cal Y}}
\newcommand{\X}{{\cal X}}
\newcommand{\T}{{\cal T}}
\renewcommand{\H}{{\cal H}}

\newcommand{\rf}[1]{(\ref{#1})}

\def\cbdu{\hfill{$\Box$}}

\newcounter{remark}
\newenvironment{remark}%
{\medskip \stepcounter{remark} \noindent {\it Remark}}{\rm \cbdu}

\title{\bf
Fractal Hamilton-Jacobi-KPZ equations}

\author{{\sc  Grzegorz Karch \& Wojbor A.~Woyczy\'nski}\\
\\
{\small Instytut Matematyczny, Uniwersytet Wroc\l awski}\\
{\small pl. Grunwaldzki 2/4, 50-384 Wroc\l aw, Poland}\\
{\small {\tt karch@math.uni.wroc.pl}}\\
{\small {\tt http://www.math.uni.wroc.pl/$\widetilde{\;}\,$karch}}\\
\\
{\small Department of Statistics and the Center for Stochastic}\\
{\small and Chaotic Processes in Science and Technology}\\
{\small Case Western Reserve University,
Cleveland, Ohio, 44106--7054, U.S.A.}\\
{\small \tt waw@po.cwru.edu} \\
{\small {\tt http://stat.cwru.edu/$\widetilde{\;}\,$Wojbor/}}\\
}

\begin{document}
\maketitle

\rightline{\it Dedicated to our friend and collaborator, Piotr
Biler}

\noindent

{\parskip = 4pt plus 4pt minus 1pt}

\begin{abstract}
Nonlinear and nonlinear evolution equations of the form $u_t=\L u
\pm|\nabla u|^q$, where  $\L$ is a pseudodifferential operator
representing the infinitesimal generator of a L\'evy stochastic
process, have been  derived as models for growing interfaces in
the case when the continuous Brownian diffusion surface transport
is augmented by  a random hopping mechanism. The goal of this
paper is to study properties of solutions to this equation
resulting from the interplay between the strengths of the
"diffusive" linear and  "hyperbolic" nonlinear terms, posed in the
whole space $\bbfR^N$, and supplemented with nonnegative, bounded,
and sufficiently regular initial conditions. 

\end{abstract}


\baselineskip=16pt

\section{Introduction}

The well-known Kardar-Parisi-Zhang  (KPZ) equation $ h_t=\nu
\Delta h +{\lambda\over 2}|\nabla h|^2 $ was derived in
\cite{KPZ86} as a model for growing random interfaces. Recall that
the interface is parameterized here by the transformation $
\Sigma(t)=(x,y, z=h(x,y,t)) $, so that $h=h(x,y,t)$ is the surface
elevation function, $\nu>0$ is identified in \cite{KPZ86} as a
``surface tension'' or ``high diffusion coefficient'', $\Delta$
and $\nabla$ stand, respectively, for the usual Laplacian and
gradient differential operators in spatial variables, and
$\lambda\in\bbfR$ scales the intensity of the ballistic rain of
particles onto the surface.

An alternative, first-principles derivation of the KPZ equation
(cf. \cite{MW01}, for more detailed information and additional
references) makes three points:

(a)   The Laplacian term can be interpreted as a result  of the
surface transport  of adsorbed particles caused by the standard
Brownian diffusion;

(b) In several experimental situations a hopping mechanism of
surface transport is present which necessitates augmentation of
the Laplacian by a nonlocal term modeled by a L\'evy stochastic
process;

(c) The quadratic nonlinearity is a result of truncation of a
series expansion of a more general, physically justified,
nonlinear even function.

\bigskip

These observations lead us to consider in this paper a nonlinear
nonlocal equation of the   form
\begin{equation}
u_t=-\L u +\lambda |\nabla u|^q, \label{eq}
\end{equation}
where
the L\'evy diffusion operator
$\L $ is defined as
\begin{equation}
\L v(x) = -\sum_{j,k=1}^N Q_{j,k} {\partial^2
v(x)\over
\partial x_j\partial x_k}
+\int_{\bbfR^N} \Big( v(x)- v(x+y)-y\cdot \nabla v(x)
1\!\!\!\;\mbox{I}_{\{|y|<1\}}(y)  \Big)\;\Pi(dy).\label{L-K-2}
\end{equation}
 The matrix
$\{Q_{j,k}\}_{j,k=1}^N$ in \rf{L-K-2} is assumed to be
a~nonnegative-definite; if it is not degenerate, a linear change
of  the variables  transforms the first term in \rf{L-K-2}
into the usual Laplacian $-\Delta$ on $\bbfR^N$ which corresponds
to the Brownian part of the diffusion modeled by $\L$. The second
term on the right-hand side of \rf{L-K-2} models the hopping
phenomena and is determined by  the Borel measure $\Pi$, usually
called the L\'evy measure of the stochastic process,  such that
$\Pi(\{0\})=0$, and $\int_{\bbfR^n}\min(1,|y|^2)\,\Pi(dy)<\infty$.
One could also include on the right-hand side a drift term $b\cdot
\nabla v$, where    $b\in\bbfR^N$ is a fixed vector but, for the
sake of the simplicity of the exposition, we omit it.  All
necessary assumptions and properties of L\'evy diffusion
operators, as well as the semigroups of linear operators generated
by $-\L$, are gathered at the beginning of  the next section.

Relaxing the assumptions that led to quadratic expression in the
classical KPZ equation, the nonlinear term in \rf{eq} has the form
$$
\lambda |\nabla u|^q=\lambda \left(|\partial_{x_1}u|^2+
...+|\partial_{x_N}u|^2\right)^{q/2},
$$
where  $q$ is a constant parameter. To study the interaction of
the "strength" of the nonlocal L\'evy diffusion parametrized by
the L\'evy measure $\Pi$, with the "strength" of the nonlinear
term, parametrized by $\lambda$ and $q$,   we consider in \rf{eq}
the whole range,   $1<q<\infty$, of the nonlinearity exponent.

Finally,  as far as  the intensity parameter $\lambda \in\bbfR$ is
concerned, we distinguish two cases:

\begin{itemize}
\item {\it The deposition case:}  Here, $\lambda>0$ characterizes
  the intensity of the ballistic deposition of particles on the evolving
interface,

\item {\it The evaporation case:} Here, $\lambda <0$, and the
model displays a time-decay of the total "mass"
$M(t)=\int_{\bbfR^N} u(x,t)\;dx$ of the solution (cf. Proposition
\ref{prop:M}).
\end{itemize}

Equation \rf{eq} will be supplemented with the nonnegative initial
datum,
\begin{equation}
u(x,0)=u_0(x), \label{ini}
\end{equation}
and our standing assumptions are that $u_0\in
W^{1,\infty}(\bbfR^N)$, and $u_0-K\in L^1(\bbfR^N)$, for some
constant $K\in \bbfR$; as usual, $W$, with some superscripts,
stands for various Sobolev spaces.

\medskip

The long-time behavior of solutions to the viscous Hamilton-Jacobi
equation $u_t=\Delta u+\lambda |\nabla u|^q$, with $\lambda\in
\bbfR$, and $q>0$, has been studied by many authors, see e.g.
\cite{AB98, BKL04,BK99, BSW02, GGK03,LS03}, and the references
therein. The dynamics of solutions to this equation is governed by
two competing effects, one resulting from the diffusive term
$\Delta u$, and the other  corresponding to the ``hyperbolic''
nonlinearity $|\nabla u|^q$. The above-cited papers aimed at
explaining  how the interplay of these two effects influences  the
large-time behavior of solutions depending on the values of $q$
and the initial data. The present  paper follows that strategy as
well.  Hence, we want to understand the interaction of the
diffusive nonlocal L\'evy  operator \rf{L-K-2} with the power-type
nonlinearity. Our results can be viewed as extensions of  some
of the above-quoted work. However,  their physical context is
quite different and, to prove them,  new mathematical tools have
to be developed.

For the sake of completeness  we  mention other
recent works  on nonlinear and nonlocal evolution equations.
First, note that equation \rf{eq} also often appears  in the context of optimal
control of
jump diffusion processes. Here, the theory of the viscosity solutions  provides
a
good framework to study these equations. We refer the reader to
the works of Jakobsen and Karlsen \cite{JK05a,JK05b}, and 
Droniou and Imbert \cite{I05,DI05} for  more
detailed information and references. Fractional conservation laws, including the
fractional Burgers equation, were studied  in
\cite{BFW98,JMW05a,JMW05b,W98} via probabilistic techniques such
as nonlinear McKean processes and interacting diffusing particle
systems.

\medskip

In the next section, we specify our assumptions on the L\'evy diffusion
operator and   state
the  main results concerning the nonlinear problem \rf{eq}-\rf{ini}. Section 3
contains proofs of those
results which are independent of the sign of the intensity parameter
$\lambda$: the existence
of solutions, the maximum principle, and the decay of $\|\nabla u(t)\|_p$
for certain $p>1$.
Further properties of solutions to \rf{eq}-\rf{ini} in the deposition case 
$\lambda>0$ are studied
in Section 4. Properties specific for the evaporation case
  $\lambda<0$  appear in Section 5. Finally, the  self-similar asymptotics of
solutions is derived
in Section 6.

\medskip

Standard notation is used throughout the paper.
For
$1\leq p\leq \infty$, the $L^p$-norm   of a Lebesgue
measurable, real-valued function $v$ defined on $\bbfR^N$
is denoted by $\|v\|_p$.
The set $C_b(\Omega)$ consists of continuous and bounded functions on $\Omega$,
and $C_b^{k}(\Omega)$ contains functions with $k$ bounded derivatives.
The space of rapidly decaying, real-valued functions is denoted by $ {\cal S}
(\bbfR^N)$.
The Fourier transform  of $v$ is 
$\widehat v(\xi)\equiv (2\pi)^{-N/2}\int_{\bbfR^N}
e^{-ix\xi} v(x)\;dx$.
The constants independent of solutions
and of $t$ (but, perhaps, dependent on the initial values) will be
denoted by the
same letter $C$, even if they may vary from line to line.
Occasionally, we write,  e.g.,  $C=C(\alpha,\ell)$ when we want to
emphasize
the dependence of $C$ on parameters $\alpha$, and $\ell$.

\section{Main results and comments}
\setcounter{equation}{0}

We begin by gathering basic properties of
solutions of the linear Cauchy problem
\begin{eqnarray}
&& u_t=-\L u ,  \quad u(x,0)=u_0(x), \label{lin:eq}
\end{eqnarray}
where the symbol $a=a(\xi)$ of the pseudodifferential operator $\L$ has the
the L\'evy--Khintchine
 representation  (cf.
\cite[Chapter 3]{J1})
\begin{equation}
 a(\xi)=Q\xi \cdot \xi
+\int_{\bbfR^n}\big(1-e^{-i\eta\xi}-{i\eta\cdot \xi}
1\!\!\!\;\mbox{I}_{\{|\eta|<1\}}(\eta)\big)\,\Pi(d\eta).\label{L-K}
\end{equation}
For every $v\in {\cal S} (\bbfR^N)$, one can use
 formula \rf{L-K} to invert the Fourier transform $\widehat {\L
 v}(\xi)=a(\xi)\widehat v(\xi)$ and to get representation \rf{L-K-2}.
In view of  \rf{L-K-2}, one can show (cf. \cite[Thm. 4.5.13]{J1}) that $\L v$
is
well defined for every $v\in C^2_b(\bbfR^N)$.

It is well-known that the operator
$-\L $
generates a positivity-preserving, symmetric  L\'evy semigroup
$e^{-t\L }$ of linear operators on $L^1(\bbfR^n)$
of the   form
\begin{equation}
(\eLt v)(x) =\int_{\bbfR^N} v(x-y)\mu^t(dy),\label{def:e}
\end{equation}
where the family $\{\mu^t\}_{t\geq 0}$ of probability Borel measures on
$\bbfR^N$ (called the convolution semigroup in
\cite{J1}) satisfies
$\widehat{\mu^t}(\xi)=(2\pi)^{-N/2}
e^{-ta(\xi)}$.
For every $1<p<\infty$, the semigroup $\eLt$ is analytic on $L^p(\bbfR^n)$, cf.
\cite[Thm. 4.2.12]{J1}. Moreover, the
representation \rf{def:e},
and the properties of the measures $\mu^t$,  imply that if $0\leq
v\leq 1$, almost
everywhere,  then  $0\leq \eLt v\leq 1$, almost everywhere (i.e., $\eLt$ is
a sub-Markovian
semigroup on $L^p(\bbfR^N)$).

The basic assumption throughout  the paper is  that the L\'evy operator
${\cal L}$ is a ''perturbation'' of the fractional Laplacian $(-\Delta)
^{\alpha/2}$, or, more precisely, that it satisfies the following  condition:
\begin{itemize}
\item
The symbol $a$ of the operator ${\cal L}$  can be written in  the form
\begin{equation}
a(\xi)=\ell|\xi|^\alpha +k(\xi),
\label{as:2}
\end{equation}
where  $\ell>0$,  $\alpha\in (0,2]$.
and  the pseudodifferential operator ${\cal K}$,
corresponding to the symbol $k$,  generates a strongly continuous
semigroup
of
operators on $L^p(\bbfR^N)$, $1\leq p\leq \infty$, with  norms uniformly bounded
in $t$.
\end{itemize}

Observe that,
 without loss of generality (rescaling   the spatial variable
$x$), we can assume that the scaling constant $\ell$ in \rf{as:2} is equal to 1.
Also, note that the above assumptions on the operator $\cal K$ 
are
satisfied if the Fourier transform of the function $e^{-tk(\xi)}$ is in
$L^1(\bbfR^N)$, for every $t>0$, and its $L^1$-norm  is uniformly bounded in
$t$.

\medskip

The  study of the large time behavior of solutions to the nonlinear
problem
\rf{eq}-\rf{ini},
will necessitate the following   supplementary asymptotic condition on $\cal L$:

\medskip

\begin{itemize}
\item The symbol
$k=k(\xi)$
appearing in  \rf{as:2}
 satisfies the condition
\begin{equation}
\lim_{\xi\to 0}{k(\xi)\over |\xi|^\alpha}=0.\label{as:3}
\end{equation}
\end{itemize}

The   assumptions \rf{as:2}
and \rf{as:3} are fulfilled,  e.g.,  by {\it
multifractional
diffusion operators\,}
$${\cal L}=-a_0\Delta+\sum_{j=1}^ka_j(-\Delta)^{\alpha_j/2},$$
with $a_0\ge 0$, $a_j>0$, $1<\alpha_j<2$, and
$\alpha=\min_{1\le j\le k}\alpha_j$, but,
more generally, one  can consider here
$$
\L=(-\Delta)^{\alpha/2}+{\cal K},
$$
where ${\cal K}$ is a generator of another L\'evy semigroup.
Nonlinear conservation laws with such nonlocal operators were studied in
\cite{BKW99,BKW01a,BKW01b}.

\medskip

In view of the assumption \rf{as:2}
imposed on its symbol $a(\xi)$,
the semigroup $\eLt$   satisfies the following decay
estimates (cf. \cite[Sec. 2]{BKW01b}, for details):
\begin{eqnarray}
\|\eLt v\|_p&\le& C
t^{-N(1-1/p)/\alpha}\|v\|_1,\label{semi}\\
\|\nabla \eLt v\|_p&\le& C
t^{-N(1-1/p)/\alpha-1/\alpha}\|v\|_1,\label{nabla}
\end{eqnarray}
for each $p\in[1,\infty]$,  all $t>0$, and a constant $C$ depending only on
$p$, and $N$. The sub-Markovian property of $\eLt$ implies
that, for every $p\in [1,\infty]$,
\begin{equation}
\|\eLt v\|_p \leq \|v\|_p.\label{contr}
\end{equation}
 Moreover, for each $p\in [1,\infty]$, we have
\begin{equation}
\|\nabla \eLt v\|_p \leq Ct^{-1/\alpha} \|v\|_p.\label{2.7a}
\end{equation}
Let us also note that  under the assumption
\rf{as:3}, the large
time behavior of $\eLt$ is described by the fundamental solution of the
linear equation
$u_t+(-\Delta)^{\alpha/2}u=0$. This results is recalled below in Lemma
\ref{lem:aslin}.

\medskip

We are now in a position to present our results concerning  the nonlinear
problem
\rf{eq}-\rf{ini}, starting with the fundamental problems of
  the existence, the uniqueness, and the regularity
of
solutions. Note that at this stage    no restrictions are imposed on the  sign
  of the parameter
$\lambda$ and  the initial datum $u_0$.
Consequently,
all results of Theorem \ref{th:exist} are valid for both the deposition, and the
evaporation cases.

\begin{theorem}\label{th:exist}
Assume that the symbol $a=a(\xi)$ of the L\'evy operator $\cal L$ satisfies
condition  \rf{as:2} with an  $\alpha\in
(1,2]$.
Then, for every $u_0\in W^{1,\infty}(\bbfR^N)$, and $\lambda\in
\bbfR$,
there exists $T=T(u_0,\lambda)$ such that
problem
\rf{eq}-\rf{ini} has a unique solution $u$ in the space
$
\X =C([0,T), W^{1,\infty}(\bbfR^N).
$

If, additionally,  there exists a constant  $K\in\bbfR$ such that $u_0-K\in
L^1(\bbfR^N)$,
 then
\begin{equation}
u-K\in C([0,T], L^1(\bbfR^N))\quad \mbox{and}\quad
\sup_{0<t\leq T}t^{1/\alpha}\|\nabla u(t)\|_1<\infty.
\label{L1-u}
\end{equation}
Moreover, for all $t\in (0,T]$,
\begin{equation}
\|u(t)\|_\infty\leq \|u_0\|_\infty\quad \mbox{and}\quad
\|\nabla u(t)\|_\infty\leq \|\nabla u_0\|_\infty.
\label{est}
\end{equation}
 and the following comparison principle is valid: for
any
two initial data satisfying condition $u_0(x)\leq \tilde u_0(x)$, the
corresponding solutions satisfy the bound $u(x,t)\leq \tilde u(x,t)$, for all
$x\in\bbfR^N$, and $t\in (0,T]$.
\end{theorem}

\bigskip

The proof of Theorem \ref{th:exist},  which is contained in Section 3, follows
the standard algorithm.
First, using the integral (mild) equation
\begin{equation}
u(t)=\eLt u_0 +\lambda \int_0^t \eLtt |\nabla u(\tau)|^q\;d\tau,
\label{duh}
\end{equation}
and the Banach fixed point
argument, we construct a  local-in-time solution.
In the next step, we prove a  ''maximum principle''  which confirms the
``parabolic
nature'' of equation \rf{eq} and allows us to prove inequalities
\rf{est}.

\bigskip

\begin{remark} {\it 2.1.}
Note that if $u$ is a solution to \rf{eq} then so is
$u-K$,  for any constant $K\in\bbfR$. Hence, without loss of generality,
in what  follows we will assume   that $K=0$.
\end{remark}

 \bigskip

\begin{remark} {\it 2.2.}
After this paper was completed  we received a preprint of \cite{DI05}
 which studied   a  nonlinear-nonlocal viscous Hamilton-Jacobi equation of the
  form
$$
u_t+(-\Delta)^{\alpha/2} u+F(t,x,u,\nabla u)=0.
$$
  Under very general assumptions on the nonlinearity, and for
$\alpha\in (1,2)$, the authors of \cite{DI05} construct a unique,
regular, global-in-time (viscosity) solution for initial data from
$W^{1,\infty}(\bbfR^N)$. Moreover, that solution also satisfies
a maximum principle which provides
inequalities \rf{est}, and the comparison principle analogous to that contained
in Theorem
\ref{th:exist}.
However, our proof of the maximum principle (cf. Theorem \ref{th:max}, below)
is simpler than the proof of the corresponding result in \cite{DI05},  and
is valid for more general L\'evy operators. On the other hand, we require the
additional assumption $u_0-K\in L^1(\bbfR^N)$, for some constant
$K\in\bbfR$.
\end{remark}

\bigskip

Once the solution $u$ is constructed, it is natural to ask questions about  its
behavior as
$t\to\infty$. From now onwards,
equation \rf{eq} will be
supplemented with the nonnegative integrable initial datum \rf{ini}.
In view of Theorem~\ref{th:exist},
the standing assumption $u_0\in W^{1,\infty}(\bbfR^N)\cap
L^1(\bbfR^N)$ allows as to define
the ``mass'' of the solution to
\rf{eq}-\rf{ini}
by the formula
\begin{equation}
M(t)=\|u(t)\|_1 =\int_{\bbfR^N} u(x,t)\;dx, \quad t\geq 0. \label{M}
\end{equation}
It's large-time behavior is one of the principal objects of study in this paper.
It turns out that in the deposition case, i.e.,  for $\lambda>0$,  the function
$M(t)$ is
increasing in $t$
(cf., Proposition \ref{prop:M}, below) and, for sufficiently small $q$,  escapes
to $+\infty$, as $t\to\infty$. More precisely, we have
the following result which is an immediate consequence of the lower bounds for
$M(t)$ obtained below in Theorem
\ref{th:M:infty:est}.

\begin{theorem}\label{th:M:infty}
Let $\lambda>0$,  $1< q\leq {N+\alpha\over N+1}$,  and
suppose that the symbol $a$ of the L\'evy operator $\cal L$ satisfies conditions
\rf{as:2} and \rf{as:3} with
$\alpha\in (1,2]$.
If  $u=u(x,t)$ is a
solution to
\rf{eq}
with an  initial datum satisfying conditions
$
0\leq u_0 \in L^1(\bbfR^N)\cap W^{1,\infty} (\bbfR^N),
$
and
$
u_0 \equiv \!\!\!\!\!\! \setminus  \; 0,
$
then
$
\lim_{t\to\infty} M(t)=+\infty.
$
\end{theorem}

When $q$ is greater that the critical exponent $(N+\alpha)/(N+1)$, we are
able to show that, for sufficiently small initial data,  the mass $M(t)$     is
 uniformly bounded in time.

\begin{theorem}\label{th:M:fin}
Let $\lambda>0$,  $ q> {N+\alpha\over N+1}$,  and
suppose that the symbol $a$ of the L\'evy operator $\cal L$ satisfies conditions
\rf{as:2} and \rf{as:3} with
$\alpha\in (1,2]$.
If,
either $\|u_0\|_1$ or $\|\nabla u_0\|_\infty$ is sufficiently small,
then
$
\lim_{t\to\infty} M(t)=M_\infty<\infty.
$
\end{theorem}

\begin{remark} {\it 2.3.}
If we limit ourselves to $\L=(-\Delta)^{\alpha/2}$ in Theorem \ref{th:M:fin}, it 
suffices only to assume that the quantity $\|u_0\|_1\|\nabla 
u_0\|_\infty^{(q(N+1)
-
\alpha-N)/(\alpha-1)}$ is small which is in perfect agreement with the 
assumption imposed in \cite{LS03} for $\alpha=2$. To see this fact, note that 
the equation $u_t=-(-\Delta)^{\alpha/2} u +\lambda |\nabla u|^q$ is invariant 
under rescaling $u_R(x,t)=R^bu(Rx,R^\alpha t)$ with $b=(\alpha-q)/(q-1)$, for 
every $R>0$. Choosing $R=\|\nabla u_0\|_\infty^{-1/(1+b)}$ we immediately obtain 
$\|\nabla u_{0,R}\|_\infty =R^{1+b}\|\nabla u_0\|_\infty=1$. Hence, the 
conclusion follows from the smallness assumption imposed on  $\| 
u_{0,R}\|_1$ in Theorem \ref{th:M:fin} and from the identity
$
\|u_{0,R}\|_1=\|u_0\|_1R^{b-n}= \|u_0\|_1\|\nabla u_0\|_\infty^{(q(N+1)-
\alpha-N)/(\alpha-1)}.
$
\end{remark}

If the L\'evy operator $\cal L$ has a non-degenerate Brownian part ,  and   if
$q\geq
2$, we
can improve Theorem \ref{th:M:fin}
showing  that the mass of every solution (not necessary small) is bounded
as $t\to\infty$.

\begin{theorem}\label{th:mass}
Let $\lambda>0$, $q\geq 2$, and suppose that the  L\'evy diffusion
operator $\cal L$
has a   non-degenerate Brownian part.  Then,  each nonnegative solution to
\rf{eq}-\rf{ini}
with an initial datum $u_0\in  W^{1,\infty}(\bbfR^N)\cap L^1(\bbfR^N)$ has
the
mass $M(t)=\int_{\bbfR^N} u(x,t)\;dx$ increasing to a finite limit
$M_\infty$,  as $t\to \infty$..
\end{theorem}

\begin{remark} {\it 2.4.}
The smallness assumption imposed in Theorem \ref{th:M:fin} seems to be
necessary. Indeed,  for $\L=-\Delta$,  it is known that if $\lambda>0$, and $
(N+2)/(N+1)<q<2 $, then
 there exists a solution to \rf{eq}-\rf{ini} such that
$\lim_{t\to\infty}M(t)=+\infty$ (cf. \cite{BSW02} and \cite[Thm. 2.4]{BKL04}).
Moreover,  if $\|u_0\|_1$ and $\|\nabla u_0\|_\infty$  are ``large'', then
the large-time behavior of solutions $u$ is dominated by the nonlinear term
(\cite{BKL04}), and  one can expect that $M_\infty=\infty$.
We conjecture that analogous results hold true at least for the $\alpha$-stable
operator (fractional Laplacian) $\L=(-\Delta)^{\alpha/2}$, and for $q$
satisfying the inequality $(N+\alpha)/(N+1)<q<\alpha $. We
also conjecture that the critical exponent $q=2$ for $\L=-\Delta$ should be
replaced
by $q=\alpha$ if $\L$ has a nontrivial $\alpha$-stable part. In this case, for
$q\geq \alpha$, we also conjecture that, as
$t\to\infty$,  the mass of any nonnegative solution   converges to a finite
limit, just like  in Theorem \ref{th:mass}. Our expectation is   that the proof
of this
conjecture can based on a reasoning similar to that contained in the proof of
Theorem \ref{th:mass}. However,  at this time, we were unable to obtain   those
estimates in a more general case.
\end{remark}
\bigskip

In the evaporation case, $\lambda <0$,  the mass $M(t)$ is a decreasing function
of
$t$ (cf.,  Proposition \ref{prop:M}, below), and the question, answered in the
next two theorems,  is when  it decays to 0 and when it decays to a positive
constant .

\begin{theorem}\label{th:M:0}
Let $\lambda<0$, $
1\leq q \leq {N+\alpha\over N+1}
$, and suppose that the symbol $a$ of the L\'evy operator $\L$ satisfies
conditions \rf{as:2} and \rf{as:3}.
If $u$ is a nonnegative solution to  \rf{eq}-\rf{ini} with an initial datum
satisfying $0\leq u_0\in W^{1,\infty}(\bbfR^N) \cap L^1(\bbfR^N)$,  then
$
\lim_{t\to\infty} M(t)=0.
$
\end{theorem}

Again, when $q$ is greater that the critical exponent, the diffusion
effects prevails for large times and,  as
$t\to\infty$,
 the mass  $M(t)$ converges to a positive limit.

\begin{theorem}\label{th:M:nonzero}
Let $\lambda<0$, $
 q> {N+\alpha\over N+1}
$, and suppose that the symbol $a$ of the L\'evy operator $\L$ satisfies
condition  \rf{as:2}.
If $u$ is a nonnegative solution to  \rf{eq}-\rf{ini} with an initial datum
satisfying $0\leq u_0\in W^{1,\infty}(\bbfR^N) \cap L^1(\bbfR^N)$,  then
$
\lim_{t\to\infty} M(t)=M_{\infty}>0.
$

\end{theorem}

The proof of Theorem \ref{th:M:nonzero} is based on the decay estimates of
$\|\nabla u(t)\|_p$ proven in Theorem \ref{grad:dec}, below. However, as was the
case for $\lambda >0$,
 we can significantly simplify   that
reasoning  for L\'evy operators
$\L$ with  nondegenerate Brownian part,  and $q\geq 2$; see the remark following
the proof of Theorem
\ref{th:M:nonzero}.

\bigskip

Our final result shows that  when the mass $M(t)$ tends  to a finite limit
$M_\infty$, as $t\to\infty$,
the solutions   to problem
\rf{eq}-\rf{ini} display a  self-similar asymptotics  dictated
by  the fundamental solution of the linear equation
$u_t+(-\Delta)^{\alpha/2}u=0$ which given by the formula
\begin{equation}
p_\alpha(x,t)=t\sp{-N/\alpha}p_\alpha (xt\sp{-1/\alpha},1)=
{1\over (2\pi)^{N/2}}
\int_{\bbfR^N} e^{ix\xi}e^{-t|\xi|^\alpha}\;d\xi.
\label{pa}
\end{equation}
More precisely, we have

\begin{theorem}\label{th:self}
Let $u=u(x,t)$ be  a  solution to problem \rf{eq}-\rf{ini}
with $u_0\in L^1(\bbfR^N)\cap W^{1,\infty}(\bbfR^N)$, and with the symbol $a$
of the L\'evy operator $\L$ satisfying conditions
\rf{as:2} and \rf{as:3}.
If
$\lim_{t\to\infty} M(t)=M_\infty$ exists and is finite, then
\begin{equation}
\lim_{{t\to\infty}}\|u(t)-M_\infty p_\alpha(t)\|_1=  0 .
\label{self-lim}
\end{equation}

If, additionally,
\begin{equation}
\|u(t)\|_p \leq Ct^{-N(1-1/p)/\alpha},
\label{self-dec-u}
\end{equation}
for some $p\in (1,\infty]$, all $t>0$, and a constant $C$ independent of $t$,
then, for every $r\in [1,p)$,
\begin{equation}
\lim_{{t\to\infty}}t^{N(1-1/r)/\alpha} \|u(t)-M_\infty p_\alpha (t)\|_r= 0 .
\label{self-lim-p}
\end{equation}

\end{theorem}

\begin{remark} {\it 2.5.} Note that, in the case $M_\infty=0$,  the results of 
Theorem \ref{th:self}
 only give that, as $t\to\infty$,  
$\|u(t)\|_r$
decays  to 0
faster than $t^{-N(1-1/r)/\alpha}$.
\end{remark}

\bigskip

\begin{remark} {\it 2.6.}
For $\lambda<0$, in view of \rf{duh}, the nonnegative solutions to
\rf{eq}-\rf{ini} satisfy 
 the
estimate $0\leq u(x,t)\leq \eLt u_0(x)$, for all $x\in\bbfR^N$, and $t>0$.
Hence, in this case,  by \rf{semi}, the decay
estimate \rf{self-dec-u} holds true with $p=+\infty$.
On the
other hand, for $\lambda>0$, the estimate of $\|\nabla u(t)\|_{p_0}$ from
Theorem \ref{grad:dec} applied to the integral equation \rf{duh} implies
immediately   \rf{self-dec-u} with $p=p_0$,  for sufficiently small
initial data; see the statement and the  proof of Theorem~\ref{grad:dec}).
In fact, following the reasoning from \cite{BK99}, it is possible to prove
\rf{self-dec-u} with $p=\infty$ without any smallness assumption. That argument
is based on the integral equation \rf{duh} and involves inequalities
\rf{semi} and \rf{est}. Here, we skip other details.
\end{remark}


\setcounter{equation}{0}
\section{Existence, uniqueness, and monotonicity}

We begin by proving that problem \rf{eq}-\rf{ini} is well-posed in
$W^{1,\infty}(\bbfR^N)$.

\begin{proposition}\label{prop:exist}
Assume that the symbol $a=a(\xi)$ of the operator $\cal L$ satisfies condition
\rf{as:2} with some $\alpha\in
(1,2]$.
Then, for every $\lambda\in \bbfR$, and $u_0\in W^{1,\infty}(\bbfR^N)$, 
there exists $T=T(u_0,\lambda)$ such that
problem
\rf{eq}-\rf{ini} has a unique solution in the space
$
\X =C([0,T), W^{1,\infty}(\bbfR^N)).
$
\end{proposition}

\proof
Our method of proof  is well-known, hence, we only sketch it.
The local-in-time solution will be constructed {\it via} the mild equation
\rf{duh}
as
the fixed point of the
operator
\begin{equation}
\T u(t)= \eLt u_0+\lambda \int_0^t \eLtt |\nabla u(\tau)|^q\;d\tau,
\label{T}
\end{equation}
in the space Banach
$\X_T =C([0,T), W^{1,\infty}(\bbfR^N))$ endowed with the norm
$$
\|u\|_{\X_T} \equiv \sup_{0\leq t\leq T} \|u(t)\|_\infty+
\sup_{0\leq t\leq T} \|\nabla u(t)\|_\infty.
$$
Inequality \rf{contr}, with $p=\infty$,  implies
\begin{eqnarray*}
\|\T u(t)-\eLt u_0\|_\infty
&\leq& |\lambda| \int_0^t \left\|\eLtt|\nabla u(\tau)\right\|^q\;d\tau\\
&\leq& C|\lambda|T\left( \sup_{0\leq \tau\leq T}
|\nabla u(\tau)\|_\infty\right)^q,
\end{eqnarray*}
and, similarly (by \rf{2.7a}, with $p=\infty$),
$$
\|\nabla \T u(t)-\nabla \eLt u_0\|_\infty \leq C|\lambda|T^{1-1/\alpha}
\left( \sup_{0\leq \tau\leq T}\|\nabla u(\tau)\|_\infty\right)^q.
$$
Moreover,    the elementary inequality
$$
|a^q-b^q|\leq C(q)|a-b|(|a|^{q-1}+|b|^{q-1})
$$
implies that, for each $R>0$, and for all $u,v\in \X_T$ such that
$\|u\|_{\X_T}\leq
R$, and
$\|v\|_{\X_T}\leq R$, we have
$$
\|\T u -\T v\|_{\X_T}\leq \left( C_1 TR^{q-1} +C_2T^{1-1/\alpha}
R^{q-1}\right)
\|u-v\|_{\X_T}.
$$
Hence, the nonlinear operator $\T$ defined in \rf{T} is a contraction on
the
ball in $\X_T$
of radius $R$ and centered at $\eLt u_0$, provided $R$ is sufficiently
large
and
$T$ is
sufficiently small.  The Banach fixed point theorem guarantees the
existence of
a solution in that ball. By a standard argument involving the Gronwall lemma,
this is the unique solution in the whole space
$\X_T$.
\cbdu

\begin{proposition}\label{prop:L1}
Assume that $u_0 \in W^{1,\infty}(\bbfR^N)$, and $u_0-K\in L^1(\bbfR^N)$,
for a constant $K\in\bbfR$.
Then the
solution
$u=u(x,t)$ of the problem \rf{eq}-\rf{ini} constructed in
Proposition~\ref{prop:exist}
satisfies condition \rf{L1-u}.
\end{proposition}

\proof
Note first that, by \rf{L-K-2},  we have $\L K=0$ for any constant
$K\in\bbfR$. Hence, replacing $u$ by $u-K$ in problem \rf{eq}-\rf{ini} one
can assume that $K=0$.
Hence,
to prove Proposition \ref{prop:L1}, it suffices to show that the operator
$\T$ used in the proof of
Proposition
\ref{prop:exist} maps the subspace of $\X_T$ defined as
$$
\Y_T\equiv \{u\in\X_T\;:\;
\sup_{0<t\leq T}\| u(t)\|_1<\infty , \quad
\sup_{0<t\leq T}t^{1/\alpha}\|\nabla u(t)\|_1<\infty\}
$$
into itself.

Observe that the properties  \rf{semi}, and \rf{nabla},
of the
L\'evy semigroup,  with $p=1$,
guarantee that $\eLt u_0\in \Y_T$,  for every $u_0\in
W^{1,\infty}(\bbfR^N)\cap L^1(\bbfR^N)$.

Now, assume that $u\in  \Y_T$.
It follows from the definition of the space $\Y_T$ that
\begin{eqnarray*}
&&\hspace{-1cm}
\left\|\int_0^t \eLtt |\nabla u(\tau)|^q\;d\tau\right\|_1\\
&\leq & \int_0^t \|\nabla u(\tau)\|_q^q\;d\tau\\
&\leq&
CT^{1-1/\alpha}
\sup_{0\leq \tau\leq T}\|u(\tau)\|_\infty^{q-1}
\sup_{0\leq \tau\leq T}\tau^{1/\alpha}\|\nabla u(\tau)
\|_1.
\end{eqnarray*}
Moreover, by \rf{nabla},  with $p=1$, we obtain
\begin{eqnarray*}
&&\hspace{-1cm}
t^{1/\alpha}\left\|\int_0^t \eLtt |\nabla u(\tau)|^q\;d\tau\right\|_1\\
&\leq &Ct^{1/\alpha} \int_0^t (t-\tau)^{-1/\alpha}
\|\nabla u(\tau)\|_q^q\;d\tau\\
&\leq&
CT^{1-1/\alpha}
\sup_{0\leq \tau\leq T}\|u(\tau)\|_\infty^{q-1}
\sup_{0\leq \tau\leq T}\tau^{1/\alpha}\|\nabla u(\tau)
\|_1.
\end{eqnarray*}
Hence $\T:\Y_T\to\Y_T$, and repeating the Banach fixed point argument in the 
space $\Y_T$, as in the proof of Proposition \ref{prop:exist}, we complete this 
proof .
\cbdu

\bigskip

Now, we are in a position to formulate and to prove the maximum principle
for
linear equations with the L\'evy diffusion  operator.
First, however, we recall an important property of positivity-preserving
semigroups and their
generators. Here, $C_0(\bbfR^N)$ denotes the space of continuous functions
decaying at infinity, and $D({\cal L})$ stands for the domain of the operator
$\cal L$.

\begin{lemma}\label{lem:positive}
Let $\L$ be a  pseudodifferential operator with the symbol $a=a(\xi)$
represented by \rf{L-K}.
Assume that $v\in D(\L)\cap C_0(\bbfR\sp N)$ and that, for some $x_0\in
\bbfR\sp N$, we have
$v(x_0)=\inf_{x\in \bbfR\sp N}v(x)\leq 0$. Then $(\L v)(x_0)\leq 0$.
\end{lemma}

\proof
This    fact is well known in the theory of generators of Feller
semigroups (cf. eg.
\cite[Ch. 4.5]{J1}) but we recall its simple proof   for the sake of
completeness of the exposition.
Since
$$
(\L v)(x_0) =\lim_{t\searrow 0} {u(x_0)-(\eLt v)(x_0) \over t},
$$
the proof will be completed by showing that $v(x_0)\leq (\eLt
v)(x_0)$,  for all $t>0$.
However, the Feller property of the semigroup $e^{-t\L}$ implies that
$\|\eLt v_{-}\|_\infty \leq \|v_{-}\|_\infty$ (cf. \rf{contr},  with
$p=\infty$). Hence
\begin{eqnarray*}
v(x_0)&=& -\|v_{-}\|_\infty
\leq -\|e^{-t\L}v_{-}\|_\infty
\leq -[e^{-t\L}v_{-}](x_0)\\
&\leq& [e^{-t\L}v_{+}](x_0)-[e^{-t\L}v_{-}](x_0)
=(e^{-t\L}v)(x_0).
\end{eqnarray*}
\cbdu

\begin{remark} {\it 3.1.}
Under the additional assumption,  $v\in D(\L)\cap C^2(\bbfR\sp N)$, it is
possible to deduce Lemma \ref{lem:positive} immediately from the
representation of the operator $\L$ given in \rf{L-K-2}. Indeed, for
$x_0\in
\bbfR\sp N$ satisfying
$v(x_0)=\inf_{x\in \bbfR\sp N}v(x)$, we have $\nabla v(x_0)=0$.
Hence,
\begin{eqnarray*}
\L v(x_0) =  -\sum_{j,k=1}^N Q_{j,k} {\partial^2
v(x_0)\over
\partial x_j\partial x_k}
+
\int_{\bbfR^N} \Big(v(x_0)- v(x_0+y)
\Big)\;\Pi(dy)\leq 0,
\nonumber
\end{eqnarray*}
since $Q$ is nonnegative-definite, $v(x_0+y)\geq v(x_0)$, for all $y\in\bbfR^N$,
and the L\'evy measure
$\Pi(dy)$ in nonnegative
(cf. also \cite[Theorem 4.5.13]{J1}).
\end{remark}

\begin{theorem}\label{th:max}
Let  $\L$ be the  L\'evy diffusion operator  defined
in \rf{L-K-2}, and
 $u=u(x,t)$, $(x,t)\in \bbfR^N\times [0,T]$,  
be a solution to  the equation
\begin{equation}
u_t+\L u+A(x,t) \cdot \nabla u=0 ,\label{eq:lin}
\end{equation}
where  $A=A(x,t)$ is a given vector field.
Moreover, suppose that the solution $u$ satisfies the following three
conditions:
\begin{equation}
u\in C_b(\bbfR^N\times [0, T]) \cap C^1_b(\bbfR^N\times [a, T]), \quad
\mbox{for each} \quad a\in (0,T),\label{as:u}
\end{equation}
\begin{equation}
u(t)\in C^2(\bbfR^N), \quad \mbox{for every} \quad t\in (0,T),
\end{equation}
and
\begin{equation}
\lim_{|x|\to\infty} u(x,t)=0,\quad \mbox{for every}\quad t\in
[0,T].\label{as:u2}
\end{equation}
Then, if  $u(x,0)\geq 0$, for all
$x\in\bbfR\sp N$,  then
$u(x,t)\geq 0$,  for all $(x,t)\in \bbfR\sp N\times [0,T]$.
\end{theorem}

\proof
For
every $t\geq 0$, 
define $f:[0,T]\to \bbfR$
by the formula
$$f(t)=\inf_{x\in\bbfR\sp N} u(x,t).$$
This is a  well-defined, and continuous function because
$u$ is uniformly
continuous and
bounded, for every $t\geq 0$ by \rf{as:u}.
Note also that, in view of  \rf{as:u2},  $f(0)=0$, and $f(t)\leq 0$,
for every $t\geq 0$.
Our goal is to show that $f\equiv 0$.

Suppose, to the contrary,  that $f(t)<0$ on an interval $(t_0,t_1)$, and
$f(t_0)=0$.
Hence, by \rf{as:u2}, for every $t\in (t_0,t_1)$, there exists an
$\xi(t)\in\bbfR\sp N$
such that
$f(t)=u(\xi(t),t)$.

Now, we show that $f(t)=u(\xi(t),t)$ is differentiable almost everywhere
on the interval
$(t_0,t_1)$,  and we follow the idea presented in \cite{CE98}.
Let us fix $s,t \in (t_0,t_1)$. If $f(t)\leq f(s)$, we obtain
\begin{eqnarray*}
0\leq f(s)-f(t)= \inf_{x\in\bbfR\sp N} u(x,s) - u(\xi(t),t)\leq
u(\xi(t),s)-u(\xi(t),t).
\end{eqnarray*}
Hence, the mean-value theorem and the assumption  \rf{as:u}
yield
$$
|f(t)-f(s)|\leq |t-s| \max_{ \tau \in [t_0,t_1]}\|u_t(\tau)\|_\infty.
$$
This means  that $f$ is locally Lipschitz on
$(t_0,t_1)$ and, 
therefore, by the Rademacher theorem,
differentiable almost everywhere . Moreover, $f'$
is bounded on every closed interval contained in $(t_0,t_1)$.

In the next step, we show that  the equality
\begin{equation}
{df(t)\over dt} =u_t(\xi(t),t)
\label{f-dif}
\end{equation}
is satisfied for   all those points from $(t_0,t_1)$ where $f$ is
differentiable.
For $t,t+h\in (t_0,t_1)$, with $h>0$, it follows from the definition of $f$
that
$$
f(t+h) =u(\xi(t+h),t+h)\leq u(\xi(t),t+h);
$$
hence
\begin{equation}
{f(t+h)-f(t)\over h} \leq {u(\xi(t),t+h)-u(\xi(t),t) \over h}. \label{f1}
\end{equation}
On the other hand,
$$
f(t-h) =u(\xi(t-h),t-h)\leq u(\xi(t),t-h),
$$
and thus, for small $h>0$,
\begin{equation}
{f(t)-f(t-h)\over h} \geq {u(\xi(t),t)-u(\xi(t),t-h) \over h}. \label{f2}
\end{equation}
 Now, we may
pass to the limit, as
$h\searrow 0$ in \rf{f1}
and \rf{f2}, to obtain the identity \rf{f-dif} in all points of
differentiability
of $f$.

To complete the proof note that under the assumption on $f$
there is a $t_2\in (t_0,t_1)$ such that
\begin{equation}
f'(t_2)<0.\label{f-es}
\end{equation}
Indeed, this follows from the fact that $0>f(t)=\int_{t_0}^t f'(s)\;ds$,  for
all $t\in (t_0,t_1)$.

However, by equality \rf{f-dif} and Lemma \ref{lem:positive}, we have
\begin{eqnarray*}
{df\over dt}(t_2)&=& -(\L u)(\xi(t_2),t_2)-a(\xi(t_2),t_2)\cdot \nabla u
(\xi(t_2),t_2)\\
&=& -(\L u)(\xi(t_2),t_2) \geq 0,
\end{eqnarray*}
which contradicts \rf{f-es}.
Hence, $f\equiv 0$ and the proof of Theorem \ref{th:max} is complete.
\cbdu

\bigskip
The maximum principle from Theorem \ref{th:max} can now be applied to our
nonlinear problem
\rf{eq}-\rf{ini}.

\begin{proposition}\label{cor:max}
The solution  $u=u(x,t)$
constructed in Theorem \ref{th:exist}
satisfies  inequalities
\rf{est},  as well as the comparison principle.
\end{proposition}

 \proof
First, we will show the comparison principle for solutions to \rf{eq}-\rf{ini} 
or, more precisely, we will prove that if 
 $u_0(x)\leq \tilde u_0(x)$, for all $x\in\bbfR^N$,  then
$u(x,t)\leq\tilde u(x,t)$, for all $(x,t)\in \bbfR^N\times [0,T]$.

Define $w=\tilde u-u$ which satisfies equation \rf{eq:lin} with
$$
A=-\lambda {|\nabla \tilde u|^q-|\nabla u|^q\over |\nabla\tilde u-\nabla
u|^2} (\nabla\tilde u-\nabla u),
$$
and with $w(x,0)=\tilde u_0(x)-u_0(x)\geq 0$.
To apply Theorem \ref{th:max}, we only need to check that $w$ satisfies
the regularity conditions imposed in \rf{as:u} - \rf{as:u2}.

Obviously, $w\in C([0,T], W^{1,\infty}(\bbfR^N))\subset C_b(\bbfR^N\times
[0,T])$. In order to improve on this statement     and  show that, actually, 
$w\in
C^1_b(\bbfR^N\times [a,T])$,  for every $a>0$, it suffices to use  the
standard bootstrap argument involving the integral equation \rf{duh}. Here, we
skip this reasoning   and   refer the reader either to
\cite[Thm.~3.1]{DI05} or to \cite[Sec. 5]{DGV03},  for more detailed 
calculations.
Next, by Theorem \ref{th:exist}, $w=\tilde u-K-u+K\in
W^{1,\infty}(\bbfR^N)\cap L^1(\bbfR^N)$ for every $t\in [0,T]$; hence,  $w$
satisfies \rf{as:u2} because
$W^{1,\infty}(\bbfR^N)\cap L^1(\bbfR^N) \subset C_0(\bbfR^N)$.

Now, the first inequality in \rf{est} follows immediately from the
comparison principle proven above because constants are solutions to
equation \rf{eq}.

To prove the second inequality in \rf{est}, we observe that the functions
$v_i=u_{x_i}$,
$i=1,...,N,$ satisfy the equations
$$
v_{i,t}=-\L v_i+\lambda q|\nabla u|^{q-2}\nabla u\cdot \nabla v_i.
$$
Applying Theorem \ref{th:max} and the reasoning from the first part of
this proof
we
obtain $$\|v_i(t)\|_\infty\leq\|v_i(0)\|_\infty,$$ which completes the
proof of Proposition
\ref{cor:max}.
\cbdu

\bigskip

\noindent {\bf Proof of Theorem \ref{th:exist}.}
The local-in-time existence of solutions is shown in Propositions
\ref{prop:exist} and \ref{prop:L1}. Proposition \ref{cor:max} provides
inequalities \rf{est} and the comparison principle.
\cbdu

\bigskip

Given $0\leq u_0\in W^{1,\infty}(\bbfR^N)\cap L^1(\bbfR^N)$, Proposition
\ref{prop:L1} allows us to  define the ``mass'' of the solution to
\rf{eq}-\rf{ini}
by the formula
$$
M(t)=\|u(t)\|_1 =\int_{\bbfR^N} u(x,t)\;dx, \quad t\geq 0.
$$
The next results shows the fundamental monotonicity property of this quantity.

\begin{proposition}\label{prop:M}
Assume that $u\in C([0, \infty), L^1(\bbfR^N))$ is a solution of problem
\rf{eq}-\rf{ini}
(or, more precisely,  a solution to the integral equation \rf{duh}).
Then, for every $t\geq 0$,
\begin{equation}
M(t)=\int_{\bbfR^N} u(x,t)\;dx =\int_{\bbfR^N} u_0(x)\;dx +\lambda
\int_0^t\int_{\bbfR^N} |\nabla u(x,\tau)|^q\; dxd\tau. \label{M2}
\end{equation}
In particular, $M(t)$ is nonincreasing   in the
evaporation
case,
$\lambda<0$, and it is nondecreasing in the deposition case,
$\lambda>0$.
\end{proposition}

\proof
Since,
for every $t\geq 0$, $\mu^t$ in the representation \rf{def:e} is a  probability
measure it follows
from the Fubini
theorem, and from the representation \rf{def:e},
that
$$
\int_{\bbfR^N} \eLt u_0(x)\; dx= \int_{\bbfR^N}\int_{\bbfR^N}
u_0(x-y)\;\mu^t(dy)dx=
\int_{\bbfR^N} u_0(y)\;dy,
$$
and, similarly,
$$
\int_{\bbfR^N} \int_0^t \eLtt|\nabla u(x,\tau)|^q\;d\tau dx =
 \int_0^t\int_{\bbfR^N}|\nabla u(x,\tau)|^q\;dxd\tau.
$$
Hence,   identity \rf{M2} is immediately obtained from equation \rf{duh} by
integrating
it
with respect
to $x$.
\cbdu

\bigskip

We conclude this section with  a result on
the time-decay 
of certain $L^p$-norms of $\nabla u$, under  smallness assumptions
on
the  initial conditions.
First, however, we need some auxiliary lemmata.

\begin{lemma}\label{lem:scat}
Let $g:(0,\infty)\to (0,\infty)$ be a  continuous function  satisfying 
the inequality $g(t)\leq A+Bg^p(t)$, for all $t>0$, some constants $A,B>0$,
and a $p>1$. If $A^{p-1}B<p^{-1}(1-1/p)^{p-1}$, and
${\rm lim\,sup}_{t\to 0}\, g(t)$ is sufficiently small, then $g(t)\leq
Ap/(p-1)$.
\end{lemma}

\proof
A direct calculation shows that the function $f(x)=x-A-Bx^p$ attains its
maximum (for $x>0$) at $x_0=(pB)^{-1/(p-1)}$. Moreover,
$f(0)=-A<0$, and
$f(x_0)=(pB)^{-1/(p-1)}(1-1/p)^{p-1}>0$,  for
$A^{p-1}B<p^{-1}(1-1/p)^{p-1}$. Hence, $g(t)$ remains in the bounded
component (containing zero) of the set
$\{x\geq 0\;:\; f(x)<0\}$. Obviously, $g(t)\leq x_0$, however, one can
improve this inequality as follows:
\begin{eqnarray*}
g(t)\leq A+Bg(t)g^{p-1}(t)\leq A+Bg(t)x_0^{p-1}\leq A+g(t)/p.
\end{eqnarray*}
This completes the proof of Lemma \ref{lem:scat}.
\cbdu

\begin{lemma}\label{lem1}
For every $p\in [1,\infty]$,
there exists a constant $C$ depending only on $p$,
 and  such that, 
for every $v\in L^p(\bbfR^N)\cap L^1(\bbfR^N)$,
\begin{equation}
D(v,p)\equiv \sup_{t>0}
t^{1/\alpha}(1+t)^\beta\|\nabla \eLt v\|_p
\leq C
\left(\|v\|_p^{1/\beta}+ \|v\|_1^{1/\beta}\right)^{\beta},
\label{Dvp}
\end{equation}
where  $\beta={N\over \alpha}\left(1-{1\over p}\right)$.
\end{lemma}

\proof
It follows from \rf{nabla} that $\|\nabla\eLt v\|_p\leq C(t/2)^{-1/\alpha}
\|e^{-(t/2)\L}v\|_p$;  hence, it suffices to estimate
$(1+t)^{N(1-1/p)/\alpha} \|\eLt v\|_p$. However, by inequalities \rf{semi},
  we have
$$
(1+t)^{N(1-1/p)/\alpha} \|\eLt v\|_p\leq \min
\left\{
(1+t)^{N(1-1/p)/\alpha}\|v\|_p,
C(p)(t^{-1}+1)^{N(1-1/p)/\alpha}\|v\|_1
\right\}.
$$
The right-hand side of the above inequality, as the function
of $t$, attains its maximum at $t_0=(C(p)\|v\|_1/\|v\|_p)^{1/\beta}$. This
completes the proof of Lemma \ref{lem1}.
\cbdu

\begin{theorem}\label{grad:dec}
Let $\lambda\in\bbfR$, $q> {N+\alpha\over N+1}$,  and 
suppose that the symbol $a$ of the L\'evy operator $\cal L$ satisfies \rf{as:2},
and
\rf{as:3},  with a certain  $\alpha\in (1,2]$.
If  $u=u(x,t)$ is a solution (not necessarily nonnegative)
to
problem
\rf{eq}-\rf{ini}, with the initial datum $u_0\in W^{1,\infty}(\bbfR^N)\cap
L^1(\bbfR^N)$, then there
  exists an  exponent $p_0$ satisfying conditions
\begin{equation}
{N+\alpha\over N+1} <p_0< {N\over N+1-\alpha},
\quad
p_0\leq q,\label{as:p}
\end{equation}
and such that,  if $D(u_0,p_0)^{p_0-1} \|\nabla u_0\|_\infty^{q-p_0}$
is sufficiently small,
then the solution $u$ satisfies the inequality
\begin{equation}
t^{1/\alpha} (1+t)^{N(1-1/p_0)/\alpha}
\|\nabla u(t)\|_{p_0}\leq C D(u_0,p_0)
\label{nab:dec}
\end{equation}
for all $t>0$, and a constant $C>0$ independent of $t$ and $u_0$.
\end{theorem}

\proof
Our reasoning is based on the integral equation \rf{duh}, estimates of the
semigroup $\eLt$ stated in \rf{semi}-\rf{2.7a},  and several algebraic
calculations on fractions. First, note that
$(N+\alpha)/(N+1)<N/(N+1-\alpha)$, for $\alpha>1$; hence the
inequalities in \rf{as:p} make sense.

In view of  equation \rf{duh} and inequalities \rf{semi}-\rf{nabla} we obtain
\begin{eqnarray}
\|\nabla u(t)\|_{p_0}&\leq& \|\nabla \eLt u_0\|_{p_0}\label{01} \\
&&+C\|\nabla u_0\|_\infty^{q-{p_0}}
\int_0^t (t-\tau)^{-N(1-1/{p_0})/\alpha -1/\alpha} \|\nabla
u(\tau)\|_{p_0}^{p_0}\;d\tau.
\nonumber
\end{eqnarray}
Next, we define the auxiliary function
$$
g(t)=\sup_{0\leq \tau\leq t} \tau^{1/\alpha} (1+\tau)^{N(1-1/p_0)/\alpha}
\|\nabla u(\tau)\|_{p_0},
$$
which, by  \rf{01},  satisfies
\begin{equation}
t^{1/\alpha}(1+t)^{N(1-1/p_0)/\alpha}\|\nabla u(t)\|_{p_0}
\leq D(u_0,p_0) +C \|\nabla u_0\|_\infty^{q-p_0} g(t)^{p_0} h(t)
\label{g:est}
\end{equation}
for all $t>0$, and a constant $C$ independent of $t$ and $u_0$. Also, let 
$$
h(t)=t^{1/\alpha} (1+t)^{N(1-1/p_0)/\alpha}
\int_0^t  (t-\tau)^{-N(1-1/p_0)/\alpha-1/\alpha}
\tau^{-p_0/\alpha} (1+\tau)^{-N(p_0-1)/\alpha}
 \;d\tau.
$$

Now, let us prove that $\sup_{t>0} h(t)<\infty$.
First,  note that, for every $t>0$,  the integral in the definition of $h(t)$
converges 
 because
the inequality $-N(1-1/p_0)/\alpha -1/\alpha>-1$ is equivalent to the condition
$p_0<N/(N+1-\alpha)$; moreover, $-p_0/\alpha>-1$   since, for    $\alpha\in
(1,2]$,
we have $p_0<N/(N+1-\alpha)\leq \alpha$.

For large values of $t$  the integral is bounded by $t^\beta$,  with
$\beta= -N(1-1/p_0)/\alpha-1/\alpha-p_0/\alpha-N(p_0-1)/\alpha+1$;  hence
$
h(t) \leq Ct^{1-p_0/\alpha-N(p_0-1)/\alpha},
$
where, for $p_0>(N+\alpha)/(N+1)$,  the exponent is negative .

Next, we analyse the behavior of $h(t)$,  as $t\to 0$. In this case, say for
$t\in (0,1)$, we obtain
$
h(t)\leq ct^{-N(1-1/p_0)/\alpha -p_0/\alpha+1}.
$
Our goal is to show that $\beta(p_0)\equiv -N(1-1/p_0)/\alpha
-p_0/\alpha+1 >0$, for each $p_0>(N+\alpha)/(N+1)$, which is sufficiently
close to $(N+\alpha)/(N+1)$. This, however, follows from the continuity of
the function $\beta(p_0)$ because,  for
$\alpha>1$, we have
$\beta((N+\alpha)/(N+1))=N(\alpha-1)^2/(\alpha(N+\alpha)(N-1))>0$.

Hence, by \rf{g:est}, the function $g(t)$ satisfies the inequality
$$
g(t)\leq D(u_0,p_0) +C \|\nabla u_0\|_\infty^{q-p_0} g(t)^{p_0},
$$
and the proof is completed by Lemma \ref{lem:scat}, because
${\rm lim\,sup}_{t\to 0}\, g(t)\leq D(u_0,p_0)$,  which follows from
\rf{g:est},  and from the properties of the function $h(t)$ shown above.
\cbdu


\section{Mass evolution in the deposition case}
\setcounter{equation}{0}

In the deposition case, i.e., for $\lambda >0$,  Proposition \ref{prop:M}
asserts that
the mass  function
$M(t)=\int_{\bbfR^N}u(x,t);dx$
is increasing in $t$.  The next results shows that, for $q\leq
(N+\alpha)/(N+1)$, as $t\to\infty$,
 the function $M(t)$ escapes to
$+\infty$
at  a certain rate, thus implying the qualitative statement of Theorem 2.2.

\begin{theorem}\label{th:M:infty:est}
Under the assumptions of Theorem \ref{th:M:infty},
there exists $T_0=t_0(u_0)$
such  that, for all $t\geq t_0(u_0)$, we have the
following lower bounds for
$M(t)$:

(a) If  $N\geq 2$, then
\begin{equation}
M(t) \geq
\left\{
\begin{array}{ccl}
C(q)\lambda M_0^q t^{(N+\alpha-(N+1)q)/\alpha}, & \mbox{for}& 1\leq q<
{N+\alpha\over
N+1};\\
C(q) \lambda M_0^q \log t, & \mbox{for}& q= {N+\alpha\over N+1}.
\end{array}
\right.
\label{+:01}
\end{equation}

(b) If $N=1$, then  
\begin{equation}
M(t) \geq
\left\{
\begin{array}{ccl}
C(q)\lambda^{1/q} M_0^{2-1/q} t^{(1+\alpha-2q)/(2q)}, & \mbox{for}& 1\leq
q<
{1+\alpha\over
2};\\
C(q) \lambda^{1/q} M_0^{q-1/q} (\log t)^{1/q}, & \mbox{for}& q=
{1+\alpha\over 2}.
\end{array}
\right.
\label{+:02}
\end{equation}
\end{theorem}

\proof
Here, we adapt the reasoning from \cite{LS03}.
Since $\lambda$ and $u_0$ are  nonnegative, it follows from equality
\rf{M2} that
$$
\lambda^{-1}M(t) =\lambda^{-1} \|u(t)\|_1\geq \int_0^t \|\nabla
u(\tau)\|_q^q\;d\tau.
$$

First,  consider $N\geq 2$.
Note that
by \rf{duh}, with $\lambda>0$,
we have
$u(x,t)\geq \eLt u_0(x)$, for all $(x,t)\in \bbfR^N\times [0,\infty)$.
Hence, by the Sobolev inequality, we obtain
\begin{equation}
\lambda^{-1}M(t) \geq C\int_0^t \|u(\tau)\|^q_{Nq/(N-q)}\;d\tau
\geq C\int_0^t \|S(\tau)u_0\|^q_{Nq/(N-q)}\;d\tau.
\label{+:M-Sob}
\end{equation}
Next, due to the assumption \rf{as:3}, we may apply Lemma \ref{lem:aslin}
from Section 6 to show
\begin{eqnarray*}
&&\hspace{-1cm}t^{N(1-1/p)/\alpha} \Big|
\|S(t)u_0\|_p-M_0\|p_\alpha(\cdot,t)\|_p\Big|\\
&&\hspace{1cm}\leq
t^{N(1-1/p)/\alpha} \|S(t)u_0-M_0p_\alpha(\cdot,t)\|_p \to 0,
\end{eqnarray*}
as $t\to\infty$. Since
$t^{N(1-1/p)/\alpha}\|p_\alpha(\cdot,t)\|_p=\|p_\alpha(\cdot,1)\|_p$
(cf. \rf{pa}),
there exists a
$t_0=t_0(u_0)$ such that
$$
\|S(\tau)u_0\|_{p}\geq {1\over 2} M_0\|p_\alpha(\cdot,1)\|_{p}
t^{-N(1-1/p)/\alpha},
\quad \mbox{for all} \quad t\geq t_0.
$$
Now, we substitute this inequality, with $p=Nq/(N-q)$, into \rf{+:M-Sob} to
obtain the
estimate
$$
\lambda^{-1} M(t)\geq C(u_0)M^q_0 \int_{t_0}^t
\tau^{-Nq(1-1/q+1/N)/\alpha}\;d\tau,
$$
which immediately implies \rf{+:01}.

The one-dimensional case requires a slightly modified argument, because
the
usual Sobolev embedding fails.
Instead, we use the interpolation inequality
$$
\|v\|^{2q-1}_\infty \leq C\|v\|_1^{q-1}\|v_x\|_q^q,
$$
for $q\geq 1$, and all $v\in L^1(\bbfR)$,  and $v_x\in L^q(\bbfR)$.
Since $\|u(t)\|_1$ is nondecreasing (cf. Proposition \ref{prop:M}), it
follows from \rf{M2}
that
\begin{eqnarray}
\|u(t)\|_1^q &\geq &
\lambda\left(\int_0^t\|u_x(\tau)\|_q^q\;d\tau \right)\|u(t)\|_1^{q-1}\\
&\geq& \lambda \int_0^t
\|u_x(\tau)\|_q^q\|u(\tau)\|_1^{q-1}\;d\tau\label{+:03}\\
&\geq & \lambda C\int_0^t\|u(\tau)\|_\infty^{2q-1}\;d\tau.\nonumber
\end{eqnarray}
Next, applying  Lemma \ref{lem:aslin} as in the case $N\geq 2$, we deduce
the
existence
of
$t_0=t_0(u_0)$, and $C>0$,  such that
$$
\|u(t)\|_\infty \geq CM_0t^{-1/\alpha}, \quad \mbox{for all} \quad t\geq
t_0.
$$
Hence, by inequality \rf{+:03}, we obtain
$$
\|u(t)\|_1^q\geq \lambda CM_0^{2q-1} \int_{t_0}^t
\tau^{-(2q-1)/\alpha}\;d\tau,
$$
which leads directly to \rf{+:02}.
\cbdu

\bigskip
At this point we are ready to provide proofs of Theorems 2.3 and 2.4.

\bigskip

\noindent{\bf Proof of Theorem \ref{th:M:fin}.}
Combining the interpolation inequality
$$
\|v\|_p\leq C(p) \|v\|_1^{N+p\over (N+1)p}
\|\nabla v\|_\infty^{N(p-1)\over (N+1)p},
$$
valid for each $p\in [1,\infty]$, and all $v\in W^{1,\infty}(\bbfR^N)\cap
L^1(\bbfR^N)$, with estimate \rf{Dvp}, we see that the quantity $D(u_0,p_0)$ 
from Theorem 
\ref{grad:dec} can 
be
controlled from above by a quantity depending only on $\|u_0\|_1$, and $\|\nabla
u_0\|_\infty$. Hence,  for small either $\|u_0\|_1$ or $\|\nabla
u_0\|_\infty$, the smallness assumption required in Theorem
\ref{grad:dec} is satisfied.

Next, the decay estimates obtained in Theorem \ref{grad:dec} allows us to prove
that $|\nabla u|^q\in L^1(\bbfR^N\times [0,\infty))$,  which immediately
implies
$M_\infty<\infty$.
Indeed, choosing $p_0$ satisfying  conditions from Theorem \ref{grad:dec},
the required integrability property of $\nabla u$ follows from the
following
inequalities
\begin{eqnarray*}
\int_0^\infty\int_{\bbfR^N} |\nabla u(x,\tau)|^q\;dxd\tau&\leq&
\|u_0\|_\infty^{q-p_0} \int_0^\infty \|\nabla
u(\tau)\|_{p_0}^{p_0}\;d\tau\\
&\leq& C \|u_0\|_\infty^{q-p_0} \int_0^\infty
\tau^{-p_0/\alpha}(1+\tau)^{-N(p_0-1)/\alpha}\;d\tau<\infty,
\end{eqnarray*}
because the condition $-p_0/\alpha-N(p_0-1)/\alpha<-1$ is automatically
satisfied   for
$p_0>(N+\alpha)/(N+1)$.
\hfill \cbdu

\bigskip

\noindent {\bf Proof of Theorem \ref{th:mass}.}
We have already mentioned in the introduction that, by a  linear change of
variables, the L\'evy operator can be written in the form
$$
\L =-\Delta +\H,
$$
where $\H$ is another L\'evy operator given  by the integral part in
\rf{L-K}.
We also recall that each L\'evy operator $\H$
is positive in the sense that, for every $p\geq 1$ and $u\in D(\H)$,  it
satisfies the inequality
\begin{equation}
\int_{\bbfR^N} (\L u) (|u|\sp{p-1}{\rm sign}\, u) \;dx \geq 0.
\label{L:poz}
\end{equation}
For the proof of \rf{L:poz}, we refer the reader to \cite[Ch. 4.6]{J1}.

In order to prove that $M_\infty<\infty$, it suffices to show that
$|\nabla u|^q
\in L^1(\bbfR^N\times [0,\infty))$. However, due to the inequality
$$
\|\nabla u(t)\|_q^q\leq \|\nabla u_0\|_\infty^{q-2} \|\nabla u(t)\|_2^2,
$$
which is  a direct consequence of \rf{est},
we only
need   to prove that
$$
\sup_{t>0} \int_0^t\int_{\bbfR^N} |\nabla u(x,\tau)|^2\;dxd\tau <\infty.
$$
For this end, we multiply equation \rf{eq} by $u^p$ and, integrating by
parts,
  obtain
\begin{eqnarray}
&&p\int_0^t\int_{\bbfR^N} u^{p-1} |\nabla u|^2\;dxd\tau +
\int_0^t \int_{\bbfR^N} u^p\H u\;dxd\tau\label{4.5.b}\\
&& ={1\over p+1} \int_{\bbfR^N} (u_0^{p+1} -u(t)^{p+1})\;dx
+\lambda \int_0^t \int_{\bbfR^N} u^p|\nabla u|^q\;dxd\tau.\nonumber
\end{eqnarray}
The second term on the left-hand side of \rf{4.5.b} is nonnegative by
inequality
\rf{L:poz}
(with  $p-1$ replaced by $p$).  Hence
\begin{equation}
\int_0^t\int_{\bbfR^N} u^{p-1} |\nabla u|^2\;dxd\tau \leq
{\|u_0\|_{p+1}^{p+1} \over p(p+1)} +{b\over p}
\int_0^t\int_{\bbfR^N} u^{p}|\nabla u|^2 \;dxd\tau,
\label{p:first}
\end{equation}
with $b=\lambda\|\nabla u_0\|_\infty^{q-2}$.

From now on, our reasoning is similar to that presented in
\cite{LS03}.
We claim that, for every integer $k\geq 1$,  
\begin{equation}
\int_0^t\int_{\bbfR^N} |\nabla u|^2 \;dxd\tau
\leq \sum_{\ell=1}^{k} {\|u_0\|_{\ell+1}^{\ell+1} b^{\ell-1}\over
(\ell+1)!}
+{b^k\over k!}
\int_0^t\int_{\bbfR^N} u^k |\nabla u|^2  \;dxd\tau.
\label{k:ind}
\end{equation}
Indeed, for $k=1$ this is just inequality \rf{p:first}, with $p=1$. To show
\rf{k:ind} for $k>1$ it is sufficient to proceed by induction.

Now, we choose $k_0$  large enough so that
$$
{b^{k_0} \|u_0\|_\infty^{k_0} \over k_0!} \leq {1\over 2}.
$$
Hence, inequality \rf{k:ind} with $k=k_0$ implies the estimate
$$
\int_0^t\int_{\bbfR^N} |\nabla u|^2 \;dxd\tau
\leq 2 \sum_{\ell=1}^{k_0} {\|u_0\|_{\ell+1}^{\ell+1} b^{\ell-1}\over
(\ell+1)!},
$$
which completes the proof of Theorem \ref{th:mass}.
\cbdu


\section{Mass evolution in the evaporation case}
\setcounter{equation}{0}

In this section, we study  equation \rf{eq} in the evaporation case, i.e., for  
$\lambda<0$.
In view of  Proposition \ref{prop:M}, the mass function $M(t)$ is now a
decreasing function
of $t$.
Our goal is  to find out under what conditions it remains bounded away from zero
or, alternatively, when  it vanishes at infinity, i.e.,  when
$M_\infty=\lim_{t\to\infty}M(t)=\lim_{t\to\infty}\int_{\bbfR^N}u(x,t)\;dx=0$.
We begin by some auxiliary results.

\begin{lemma}\label{lem:1-}
If   $w\in W^{1,1} (\bbfR^N)$ then, for every $R>0$,
$$
\|w\|_1\leq 2R \int_{|x|\leq 3R} |\nabla w(x)|\;dx +2
\int_{|x|>R}|w(x)|\;dx.
$$
\end{lemma}

A short and elementary proof of Lemma \ref{lem:1-} can be found in the
paper by
Ben-Artzi
and Koch \cite{BK99}.

\begin{lemma}\label{lem:2-}
Let $\lambda <0$. Assume that the symbol $a=a(\xi)$ of the L\'evy operator $\cal
L$ given by the formula \rf{L-K} satisfies
the assumptions \rf{as:2} and \rf{as:3}.
If  $r\in C(0,\infty)$ is  a nonnegative function such that
\begin{equation}
\lim_{t\to\infty} r(t)t^{-1/\alpha}=+\infty,
\label{r:as}
\end{equation}
then
$$
\lim_{t\to\infty}\int_{|x|\geq r(t)} u(x,t)\;dx =0.
$$
\end{lemma}

\proof
Since $\lambda <0$, it follows from the integral equation \rf{duh}   that,
for all $t\geq 0$,  and $x\in\bbfR^N$, we have
$
0\leq u(x,t)\leq \eLt u_0(x)
$. Hence
\begin{eqnarray*}
\int_{|x|\geq r(t)} u(x,t)\;dx &\leq & \int_{|x|\geq r(t)} \eLt
u_0(x)\;dx\\
&\leq & \int_{|x|\geq r(t)} \Big| \eLt u_0(x)-\|u_0\|_1
p_\alpha(x,t)\Big|\;dx\\
&&+
\|u_0\|_1 \int_{|x|\geq r(t)} p_\alpha(x,t)  \;dx.
\end{eqnarray*}
As $t\to\infty$, the first term on the right-hand side tends to 0  by Lemma
\ref{lem:aslin}  with
$p=1$. In the second term, we change the variables $y=xt^{-1/\alpha}$ to
obtain
$$
\int_{|x|\geq r(t)} p_\alpha(x,t)\;dx = \int_{|y|\geq r(t)t^{-1/\alpha}}
p_\alpha(y,1)\;dy \to 0,
\quad\mbox{as}\quad t\to\infty,
$$
in view of the self-similarity of the form $p_\alpha(x,t)=t^{-N/\alpha}p_\alpha
(xt^{-1/\alpha},1)$, the assumption \rf{r:as},  and  since $p_\alpha
(\cdot,1)\in L^1(\bbfR^N)$.
\cbdu

\bigskip

Now we are ready to prove the results of Section 2 describing   mass evolution
in  the evaporation case.

\bigskip

\noindent {\bf Proof of Theorem \ref{th:M:0}.}
Since $M(t)$ is nonnegative,  equation
\rf{M2} with $\lambda<0$ implies that
$$
\int_0^\infty \int_{\bbfR^N} |\nabla u(x,\tau)|^q\;dxd\tau \leq \|u_0\|_1.
$$
For $t\geq 0$, define
$$
\omega (t)=\left(\int_{t/2}^\infty \int_{\bbfR^N} |\nabla
u(x,\tau)|^q\;dxd\tau
\right)^{1/q},
$$
 and   notice that $\omega=\omega(t)$ is a nonincreasing
function
on
$[0,\infty)$ which satisfies condition
\begin{equation}
\lim_{t\to\infty} \omega(t) =0. \label{omega:decay}
\end{equation}
Now, for $t\geq 1$, $s\in (t/2,t)$, and  $R>0$, we infer from Lemma
\ref{lem:1-}
combined with
the H\"older inequality that
\begin{equation}
\|u(s)\|_1\leq CR^{1+N(1-1/q)}\|\nabla u(s)\|_q +2 \int_{|x|>R}
|u(x,s)|\;dx.\label{-:01}
\end{equation}
Since $s\mapsto \|u(s)\|_1$ is  nonincreasing on $(t/2,t)$, it
follows from inequality
\rf{-:01}, and the H\"older inequality, that
\begin{eqnarray}
\|u(t)\|_1&\leq& {2\over t} \int_{t/2}^t \|u(s)\|_1\;ds\label{-:02}\\
&\leq & CR^{1+N(1-1/q)} t^{-1/q} \omega (t)
+{4\over t} \int_{t/2}^t \int_{|x|>R} |u(x,s)|\;dxds.\nonumber
\end{eqnarray}
Next,   fix  $\delta \in (0, (1+N(1-1/q))^{-1})$ and 
define
$$
R(t)=t^{1/\alpha} \omega (t)^{-\delta}.
$$
This function is nondecreasing
(because $\omega(t)$ is nonincreasing) which implies that $R(t)\geq
R(s)$
for all
$s\in [t/2,t]$. Hence, substituting $R=R(t)$ into inequality \rf{-:02} we
obtain
\begin{eqnarray*}
\|u(t)\|_1&\leq& Ct^{(1+N(1-1/q))/\alpha -1/q}
\omega(t)^{1-\delta(1+N(1-1/q))}\\
&&+
{4\over t} \int_{t/2}^t \int_{|x|>R(s)} |u(x,s)|\;dxds.
\end{eqnarray*}
The first term on the right-hand side tends to zero, as $t\to\infty$,
because the
inequality
$(1+N(1-1/q))/\alpha -1/q\leq 0$ is equivalent to $q\leq (N+\alpha)/(N+1)$
and
$\omega(t)^{1-\delta(1+N(1-1/q))}\to 0$. The second term converges to zero
by the Lebesgue
Dominated
Convergence
theorem,  and Lemma \ref{lem:2-},  because
$$
{4\over t} \int_{t/2}^t \int_{|x|>R} |u(x,s)|\;dxds=
{4} \int_{1/2}^1 \int_{|x|>R(t\tau)} |u(x,t\tau)|\;dxd\tau.
$$
Hence $M(t)=\|u(t)\|_1\to 0$, as $t\to\infty$, and the proof is complete.
\cbdu

\bigskip

\noindent{\bf Proof of Theorem \ref{th:M:nonzero}.} For $\e\in (0,1]$, we
denote by $u^\e=u^\e(x,t)$ the solution to \rf{eq}-\rf{ini} with $\e
u_0(x)$ as the initial datum. Since, by the comparison principle from
Theorem \ref{th:exist},
$0\leq u^\e (x,t)\leq u(x,t)$, for all $x\in\bbfR^N$, and $t>0$, it suffices
to show that, for small $\e>0$,
$$
\lim_{t\to\infty} M^\e(t)=\lim_{t\to\infty}\int_{\bbfR^N}
u^\e(x,t)\;dx=M^\e_\infty>0.
$$
However, by Proposition \ref{prop:M},
$$
M^\e_\infty= \e \int_{\bbfR^N} u_0(x)\;dx -|\lambda| \int_0^\infty
\int_{\bbfR^N} |\nabla u^\e(x,\tau)|^q\;d\tau.
$$
Hence, for sufficiently small $\e>0$, and for $p_0$ satisfying the
assumptions of Theorem \ref{grad:dec}, we have
\begin{eqnarray*}
M^\e_\infty
&\geq&  \e \int_{\bbfR^N} u_0(x)\;dx \\
&&-|\lambda| \|\nabla\e u_0\|_\infty^{q-p_0}
\int_0^{\infty} \|\nabla u^\e(\tau)\|_{p_0}^{p_0}\;d\tau\\
&\geq&  \e \int_{\bbfR^N} u_0(x)\;dx \\
&&-|\lambda| C\|\nabla\e u_0\|_\infty^{q-p_0}D(\e
u_0,p_0)^{p_0}\int_0^{\infty}\tau^{-p_0/\alpha}
(1+\tau)^{-N(p_0-1)/\alpha}\;d\tau.
\end{eqnarray*}
We have shown already in the proof of Theorem \ref{grad:dec} that the
integral on the right-hand side is finite. Moreover, by the definition of
$D(u_0,p_0)$, we have $D(\e u_0,p_0)=\e D(u_0,p_0)$ and, consequently,
\begin{equation}
M_\infty^\e \geq \e \int_{\bbfR^N} u_0(x)-\e^q|\lambda|
\|\nabla u_0\|_\infty^{q-p_0} D(u_0,p_0)C,\label{*}
\end{equation}
with $C>0$ independent  of $\e$,  and $u_0$. Since $q>1$, it follows
from \rf{*} that,  for sufficiently small $\e>0$,  necessarily  $M_\infty^\e>0$.
\cbdu

\bigskip

\begin{remark} {\it 5.1.} If we strengthen the assumptions
   in Theorem \ref{th:M:nonzero},  and demand that  $q\geq 2$,  and 
$\L$ has a  nondegenerate Brownian part,  {i.e.},  $\L =\-\Delta +{\cal H}$,
(cf.,  the proof of Theorem \ref{th:mass}), then,
multiplying equation \rf{eq} by $u$,  and integrating over $\bbfR^N\times
[0,t]$,  we obtain
\begin{eqnarray*}
&&2\int_0^t\int_{\bbfR^N}  |\nabla u|^2\;dxd\tau +
\int_0^t \int_{\bbfR^N} u\H u\;dxd\tau\\
&& ={1\over 2} \int_{\bbfR^N} (u_0^{2} -u(t)^{2})\;dx
-|\lambda|\int_0^t \int_{\bbfR^N} u|\nabla u|^q\;dxd\tau.\nonumber
\end{eqnarray*}
In particular, we have (cf. \rf{L:poz} and \rf{4.5.b})
$$
\int_0^t\int_{\bbfR^N}  |\nabla u|^2\;dxd\tau\leq {1\over 2} \|u_0\|_2^2.
$$
Hence, repeating the reasoning from the proof of Theorem
\ref{th:M:nonzero}, we obtain, for sufficiently small $\e>0$,  that
$$
M_\infty^\e \geq \e \int_{\bbfR^N} u_0(x)\;dx -\e^q (|\lambda|/2) \|\nabla
u_0\|_\infty^{q-2}\|u_0\|_2^2>0.
$$
Note that, in this case,  we do not need
decay estimates from Theorem \ref{grad:dec}.
\end{remark}

\section{Self-similar asymptotics}
\setcounter{equation}{0}

Assumptions   \rf{as:2} and \rf{as:3} allow us  to approximate $\eLt u_0$
by a
multiplicity
of the kernel
$$
p_\alpha(x,t)={1\over (2\pi)^{-N/2}}
\int_{\bbfR^N} e^{ix\xi}e^{-t|\xi|^\alpha}\;d\xi.
$$
Indeed, we have the following

\begin{lemma}\label{lem:aslin}
If the symbol $a(\xi)$ of the L\'evy operator $\cal L$ satisfies conditions
\rf{as:2},  and \rf{as:3}, 
then, for   each $p\in [1,\infty]$, and $u_0\in L^1(\bbfR^n)$,
\begin{equation}
\lim_{t\to\infty} t^{N(1-1/p)/\alpha}
\left\|e^{-t\cal L}u_0- \int_{\bbfR^N}u_0(x)\;dx
\;p_\alpha(t)\right\|_p=0.
\label{lim:lin1}
\end{equation}
\end{lemma}

\proof
This result is obtained immediately  from the inequality
\begin{equation}
\left\|h\ast g(\cdot)-\left(\int_{\bbfR^n} h(x)\,
dx\right)g(\cdot)\right\|_p
\le C\|\nabla g\|_p\|h\|_{L^1(\bbfR^n,|x|\, dx)}, \label{ineq}
\end{equation}
which is valid for each $p\in[1,\infty]$, all $h\in L^1(\bbfR^n,|x|\, dx)$, and
every
$ g\in C^1(\bbfR^n)\cap W^{1,1}(\bbfR^n)$, with
a constant $C=C_p$ independent of $g$,  and $h$.
The  inequality itself is a simple consequence of the Taylor expansion.

To prove the Lemma  we apply \rf{ineq}, with
$h=e^{-t{\cal K}}u_0$, and $g(x)=p_\alpha(x,t)$,
assuming first that $u_0\in L^1(\bbfR^n,|x|\, dx)$. The general case
of
$u_0\in L^1(\bbfR^n)$ can then be handled by an approximation argument.
Details of such a reasoning can be found  in \cite[Cor. 2.1 and 2.2]{BKW01b}.
\cbdu

\bigskip

Now, we are in a position to prove the final, main theorem of Section 2.

\bigskip

\noindent{\bf Proof of Theorem \ref{th:self}.}
Since $M_\infty$ is finite, formula \rf{M2} for $M(t)$ implies that
$\|\nabla
u\|_q^q\in L^1(0,\infty)$.
It follows from the integral equation \rf{duh} that
$$
\|u(t)-e^{-(t-t_0)\L}u(t_0)\|_1\leq \int_{t_0}^t\|\nabla
u(\tau)\|_q^q\;d\tau,
\quad \mbox{for all} \quad t\geq t_0\geq 0.
$$
Hence,
\begin{eqnarray*}
\|u(t)-M_\infty p_\alpha(t)\|_1&\leq& \int_{t_0}^t \|\nabla
u(\tau)\|_q^q\;d\tau
\\
&&+\|e^{-(t-t_0)\L}u(t_0)-M(t_0)p_\alpha(t)\|_1\\
&&+
\|p_\alpha(t)\|_1\Big| M(t_0)-M_\infty\Big|.
\end{eqnarray*}
Letting $t\to\infty$, and using Lemma \ref{lem:aslin}
for  the second term on the
right-hand side,  we obtain
$$
\limsup_{t\to\infty} \|u(t)-M_\infty p_\alpha(t)\|_1 \leq
\int_{t_0}^\infty
\|\nabla u(\tau)\|_q^q\;d\tau
+|M(t_0)-M_\infty|,
$$
because $\|p_\alpha(t)\|_1=1$. Since $t_0$ can be arbitrarily large, and
each term
on the right-hand side
tends to 0, as $t_0\to\infty$,   the proof of \rf{self-lim} is complete.

Now, we apply the
H\"older inequality  to obtain
$$
\|u(t)-M_\infty p_\alpha (t)\|_r \leq
\|u(t)-M_\infty p_\alpha (t)\|_1^{1-\gamma} \left(
\|u(t)\|_p^{\gamma}+M_\infty\|p_\alpha (t)\|_p^{\gamma}\right),
$$
with $\gamma =(1-1/r)/(1-1/p)$.
Finally, to prove \rf{self-lim-p},
it suffices to apply inequality \rf{self-dec-u}, the asymptotic result
\rf{self-lim},  and the identity
$\|p_\alpha(t)\|_p=t^{-N(1-1/p)/\alpha} \|p_\alpha (1)\|_p$.
\cbdu

\bigskip

{\bf Acknowledgements.}
The authors would like to thank J. Droniou for making his unpublished manuscript
  \cite{DI05} available to them.
This paper was partially written  while the first-named
author enjoyed the hospitality and support 
of the Center for
Stochastic and Chaotic Processes in Science and Technology at Case Western
Reserve University, Cleveland, Ohio, U.S.A., sponsored by the U.S. National
Science Foundation Grant  INT-0310055, and
of the Helsinki University
of Technology,  and the University of Helsinki, Finland, within the
Finnish Mathematical Society Visitor Program in Mathematics 2005-2006,
{\sl Function Spaces and Differential Equations}.
The preparation of this paper was also partially supported by the KBN grant
2/P03A/002/24, 
and by the European Commission Marie Curie Host Fellowship
for the Transfer of Knowledge ``Harmonic Analysis, Nonlinear
Analysis and Probability''  MTKD-CT-2004-013389.

\medskip


\end{document}